\documentclass[graybox]{svmult}

\usepackage{type1cm}        
\usepackage{graphicx}        
\usepackage[bottom]{footmisc}

\usepackage{newtxtext}       %
\usepackage{newtxmath}       

\usepackage[utf8]{inputenc}
\usepackage[T1]{fontenc}
\usepackage{amsmath}
\usepackage{hyperref}
\usepackage{svg}
\usepackage{placeins}
\usepackage[backend=biber]{biblatex}
\graphicspath{{figures/}}
\usepackage{mathtools}
\usepackage{listings}
\usepackage{setspace}
\usepackage{multirow}
\usepackage{longtable}
\usepackage{floatrow}
\begin{document}

\title*{Surrogate Models for Optimization of Dynamical Systems
\thanks{This research was supported by Joint MSc in Applied and Interdisciplinary Mathematics,
coordinated by the University of L’Aquila (UAQ) in Italy, Department of Information Engineering, Computer Science and Mathematics (DISIM) and the Deutsche Forschungsgemeinschaft through the International Research Training Group 1792 "High Dimensional Nonstationary Time Series", the Yushan Scholar Program of Taiwan, and European Union’s Horizon 2020 training and innovation programme ”FIN-TECH”, under the grant No. 825215 (Topic ICT-35-2018, Type of actions: CSA).}}
\author{Kainat Khowaja, Mykhaylo Shcherbatyy, and Wolfgang Karl Härdle}
\institute{Kainat Khowaja \at  International Research Training Group 1792 "High Dimensional Nonstationary Time Series", Humboldt-Universität zu Berlin, Berlin, Germany; Ivan Franko National University of Lviv, Ukraine; University of L'Aquila, Italy. \email{kainat.khowaja@hu-berlin.de}
\and Mykhaylo Shcherbatyy \at Ivan Franko National University of Lviv, Ukraine. \email{mykhaylo.shcherbatyy@lnu.edu.ua}
\and Wolfgang Karl Härdle\at BRC Blockchain Research Center, Humboldt-Universität zu Berlin,  Germany; Sim Kee Boon Institute, Singapore Management University, Singapore; WISE Wang Yanan Institute for Studies in Economics, Xiamen University, China; Dept. Information Science and Finance, National Chiao Tung University, Taiwan, ROC; Dept. Mathematics and Physics, Charles University, Czech Republic. Grants–DFG IRTG 1792, CAS: XDA 23020303  and COST Action CA19130 gratefully acknowledged. \email{haerdle@hu-berlin.de}}
%
%
\maketitle
\let\cleardoublepage\clearpage 
\vspace*{-25mm}
\abstract{
Driven by increased complexity of dynamical systems, the solution of system of differential equations through numerical simulation in optimization problems has become computationally expensive. This paper provides a smart data driven mechanism to construct low dimensional surrogate models. These surrogate models reduce the computational time for solution of the complex optimization problems by using training instances derived from the evaluations of the true objective functions. The surrogate models are constructed using combination of proper orthogonal decomposition and radial basis functions and provides system responses by simple matrix multiplication. Using relative maximum absolute error as the measure of accuracy of approximation, it is shown surrogate models with latin hypercube sampling and spline radial basis functions dominate variable order methods in computational time of optimization, while preserving the accuracy. These surrogate models also show robustness in presence of model non-linearities. Therefore, these computational efficient predictive surrogate models are applicable in various fields, specifically to solve inverse problems and optimal control problems, some examples of which are demonstrated in this paper.}


\chapter{Introduction}
Over the years, mathematical modeling and optimization techniques have effectively described complex real-life dynamical structures using system of differential equations. More often, the dynamical behavior of such models, especially in optimization and inverse problems (the problems where some of the 'effects' (responses) are known but not some of the 'causes' (parameters) leading to them are unknown), cause necessity of repetitive solution of these model equations with a slight change in system parameters. While numerical models replaced experimental methods due to their robustness, accuracy, and rapidness, their increasing complexity, high cost, and long simulation time have limited their application in domains where multiple evaluations of the model differential equations are demanded. 

To prevent this trade-off between computational cost and accuracy, one needs to focus on Reduced Order Models (ROM) which provide compact, accurate and computationally efficient representations of ODEs and PDEs to solve these multi-query problems. These approximation models, also commonly recognized as a surrogate models or meta-models \cite{shcherbatyy2018}, allow the determination of solution of model equations for any arbitrary combination of input parameters at a cost that is independent of the dimension of the original problem. Accordingly, they meet the most essential criteria of every analysis problem:  the criteria of highest fidelity at lowest possible computational cost, where high fidelity is defined by the efficacy of theoretical methods to replicate the physical phenomenons with least possible error \cite{Emiliano2013}.

This paper employs Proper Orthogonal Decomposition (POD), a model reduction technique originating in statistical analysis and known for its optimality as it captures the most dominant components of data in the most efficient way \cite{hinze}. POD serves the purpose of dimension reduction by extracting hidden structures from high dimensional data and projecting it on lower dimensional space \cite{springer2005}. In this work, POD will be used to derive low order models of dynamical system by reducing a high number of interdependent variables to a much smaller number of uncorrelated variables while retaining as much as possible of the variation in the original variables. 

Over a century ago, Pearson proposed the idea of representing the statistical data in high dimensional space using a straight line or plane, hence discovering a finite dimensional equivalence of POD as a tool for graphical analysis \cite{Pearson1901}. In the years following Pearson's paper, the technique has been independently rediscovered by several other scientists including Kosambi, Hotelling and Van Loan under different names in the literature such as Principle Component Analysis (PCA), Hotelling Transformation and Loeve-Karhunen Expansion, depending on the branch in which it is being tackled. Despite its early discovery, the availability of computational resources required to compute POD modes were limited in earlier  years and the technique remained virtually unused until 1950s. The technological advancements took an upturn after that with the invention of powerful computers and led to the popularity of POD \cite{springer2005}. Since then, the development and applications of POD have been widely investigated in diverse disciplines such as structural mechanics \cite{springer2005}, aerodynamics \cite{Emiliano2013}, signal and image processing \cite{Benaarbia2017}, etc. Due to its strong theoretical foundations, the technique has been used in many applications, such as for damage detection \cite{Lanata2006}, human face recognition \cite{Kirby1987}, detection of signals in multi-channel time-series \cite{Wax1985}, exploration of peak clustering \cite{berardi2015} and many more.

In general, a non-equivalent variant of POD, known as factor analysis, has been renowned and has been used for various applications \cite{Felix2018,Anita2015, Bai2016, Alessandro2018}, etc. Unlike POD, factor analysis assumes that the data have a strict factor structure and it looks for the factors that amount for common variance in the data. On contrary, PCA the finite counterpart of POD, allows the accountability of maximal amount of variance for observed variables. The PCA analysis consists of identifying the set of variables, also known as principle components, from the system that retain as much variation from the original set of variables as possible. Similarly, Principal Expectile Analysis (PEC), which generalizes PCA for expectiles was recently developed as a dimension reduction tool for extreme value theory \cite{Haerdle2019}. These POD equivalent tools have also been adopted in analysis on several instances such as \cite{Felix2018, Ying2018, Li2018, Haerdle2019}. Yet, most of the literature exploits only the real life data for dimension reduction. Even though some analysis highly relies on real life data, there is an urgent need of introduction of tools that utilize simulated data generated from the non-standard models with nonlinear differential equations that are on constant rise and hold potential for enrichment of analysis.

Moreover,  optimal control problems and mathematical optimization models are widely seen in various applications. These models are often used for normative purposes to solve the minimization/maximization problems and require repetitive evaluation in various context with different parameter values to find the optimum set of parameters. This parameter exploration process can be computationally intense, specially in complex non-linear system which emphasize the need of dimension reduction for these models. 

Through this research, the application of POD to reduce the dimensionality of dynamical systems is proposed. The present work resorts to explore the efficacy of POD on few common applications, the models which have been previously defined and commonly used. We hypothesize that the system responses of dynamical models can be obtained with a very high accuracy, but lower computational cost model reduction technique. The novelty of this hypothesis lies in the fact that dimensional reduction techniques have rarely been explored for optimal control problems, specially the combination of POD and Radial Basis Functions (RBF) to make surrogate models is quite under utilized, specially for the the models discussed in this paper. 

The computational procedure of the research is decomposed between offline and online phases. The offline phase (training of the model) entails utilization of sampling techniques to generate data, computation of snapshot matrix of  model solutions using variable order methods for solving of ODE (model of dynamical system), obtainment of proper orthogonal modes via Singular Value Decomposition (SVD) and estimation of POD expansion coefficients that approximate the POD basis (via interpolation techniques radial basis functions). The online phase (testing of the model) involves redefinition of model equations in terms of surrogate models and computation of system responses corresponding to any arbitrary set of input parameters in given domain \cite{shcherbatyy2018}. Next, the quality of the model is validated and evaluated by carrying out error analysis and various experimental designs are employed by varying sampling and interpolation techniques and changing the size of training set to determine the combination that generates that results in the least maximum absolute error. Finally, using that experimental design, optimization criterion are calculated using both models to evaluate accuracy of the model. For the computations, a MATLAB software is developed by the author which utilizes a combination of inbuilt and user-defined functions. The illustrations used in this work are also generated using MATLAB.

The next chapter will lay down theoretical concepts related to POD, SVD and RBF, and how surrogate models are constructed to project the dynamical system onto subspaces consisting of basis elements that contain characteristics of the expected solution. Chapter 2  will explain how the computational procedure (algorithm) and Chapter 3 will implement the concepts developed in Chapter 1 and methodology presented in Chapter 2 on a set of dynamical systems. Finally, last chapter will conclude the main results and provide a summary of current research, its limitations, as well as future prospects.
\chapter{Mathematical Formulation}
\label{MathematicalFormulation}
\setcounter{section}{0}

Model reduction techniques have been known for their ability to reduce the computational complexity of mathematical models in numerical simulations. The main reason ROM has found applications in various disciplines is due to its strong theoretical foundations and the demand of model reduction techniques in ever-so-rising computational complexities and intrinsic property of high dimensionality of physical system. ROM addresses these issues effectively by providing low dimensional approximations.

Although a variety of dimensionality-reduction techniques exist such as operational based reduction methods \cite{Schilders2008}, reduced basis methods  \cite{Boyaval2010}, the ROM methodology is often based upon POD. Analogous to PCA, the POD theory requires to find components of the systems, known as Proper Orthogonal Modes (POMs), that are ordered in a way that each subsequent mode holds less energy than previous one. As stated earlier, POD is ubiquitous in the dimensionality reduction of physical systems. It presents the optimal technique for capturing the system modes in least square sense. That is, for constructing ROM for any system, incorporating \textit{k} POMs will give the best \textit{k} component approximation of that system. This assures that any approximation formulated using POD will be the best possible approximation: there is no other method that can reduce the dimensionality of the given system in lower number of components or modes. 

This chapter discusses in depth the mathematical concepts associated with POD and its correspondence with SVD and RBF for construction of surrogate models. The computational procedure adapted in Chapter \ref{chapter-computation} and Chapter \ref{chapter-applications} is strictly based on the theory formulated in this chapter.

\section{Formulation of Optimization Problem}
Many problems of optimal control are focused on the minimization and maximization problems. In order to find an optimal set of parameters, optimization models are usually defined in which the problems are summarized by the objective function. These optimization parameters are called control parameters and they affect the choice of allocation. In optimal control problems, these parameters are time paths which are chosen within certain constraints so as to minimize or maximize the objective functional. The applications presented in Chapter \ref{chapter-applications} are optimization problems, the general structure of which has been discussed in the next paragraph.

Let us consider optimization problem which consists of finding a vector of optimization parameters $u^* \in U_{S}$ and proper state function $y^* \subset Y_{S}$, that minimizes the optimization criterion (objective function) 

\begin{equation}
\psi_0 = \tilde\psi_0(u^*,y^*) = \min_{(u,y)\in U_S \times Y_S} \tilde\psi_0(u,y)
\label{optimizationequations}
\end{equation}

subject to ODEs (state equation)
 
\begin{equation}
c(y,u)= 0 \sim \begin{dcases}
y_i'-f(t,u,y)=0,  \ t\in [t_0,T],\\
y(t_0)-y_0=0,
\end{dcases}
\label{stateequations}
\end{equation}

box constrains on the control variable

\begin{equation}
    U = \{u \in U_S: u^- \leq u \leq u^+, u^-\in U_S, u^+\in U_S\}
\end{equation}

and possibly additional equality and non-equality constraints on state and control

\begin{equation}
\begin{matrix}
 \tilde{\psi_j}(u,y)=0, j= 1, \ldots,m_1,\\
\tilde{\psi_j}(u,y)\leq 0, j= m_1+1, \ldots,m.
\end{matrix}
\label{eq:constraints}
\end{equation}

where $U_S$ and $Y_S$ are real Banach spaces, $u= u(t)= [u_1(t), \ldots, u_{n_u}(t)]^\top \in U_S, y=y(t)= [y_1(t), \ldots,  y_{n_y}(t)]^\top \in Y_S, \tilde{\psi_j}:U_S\times Y_S \rightarrow \mathbb{R}, j=0,1,\ldots,m$

We assume that for each $u\in U$, there exists a unique solution $y(u)$ of state equation $c(y,u)=0$. With the aim of compactness, we will write optimization problem (\ref{optimizationequations}- \ref{eq:constraints}) in reduced form: find a function $u^*$ such that

\begin{equation}
\begin{matrix}
u^{*} \in U_{\partial_u}, \psi_{0}\left(u_{*}\right)=\displaystyle \min_{u \in U_{\partial_u}} \psi_{0}(u)\\
U_{\partial_u}=\left\{u: u \in U ; \psi_{j}(u)=0, j=1, \ldots, m_{1} ; \psi_{j}(u) \leq 0, j=m_{1}+1, \ldots, m\right\}\\
c(y(u), u)=0\\
\psi_{j}(u)=\tilde{\psi}_{j}(u, y(u)), j=0,1, \ldots, m
\end{matrix}
\end{equation}
The optimal control problems in this research are solved using direct method. Each problem is transformed to nonlinear programming problem, i.e., it is first discretized and then the resulting nonlinear programming problem is optimized. The advantage of direct methods is that the optimality conditions of an non linear programming problems are generic, whereas optimality conditions of undiscretized optimal control problems need to be reestablished for each new problem and often require partial a-priori knowledge of the mathematical structure of the solution which in general is not available for many practical problems.

The first step in the direct method is to approximate each component of the control vector by a function of finite parameters $u_i(t)= u_i(t,b^{(i)}), b^{(i)} = [b^{(i)}_1,...,b^{(i)}_{n_i}]^\top, i= 1, \ldots, n_u$. As a result, we write control function $u(t)$ as a function of vector of optimization parameters $b$: $u(t)=u(t,b)$. In this paper we use a piecewise-linear or piecewise-constant approximation for each function $u_{i}(t), i=1, \ldots, n_{u}$.

The optimization problem can be written as nonlinear programming problem in such a way that we have to find a vector $b^*$ such that

\begin{equation}
\begin{matrix}
b^{*} \in U_{\partial}, \psi_{0}\left(b^{*}\right)= \displaystyle \min _{b \in U_{\partial}} \psi_{0}(b) \\
U_{\partial}=\left\{b: b \in U_{b}, \psi_{j}(b)=0, j=1, \ldots, m_{1} ; \psi_{j}(b) \leq 0, j=m_{1}+1, \ldots, m\right\} \\
U_{b}=\left\{b: b \in R^{n}, b^{-} \leq b \leq b^{+}, b^{-} \in \mathbb{R}^{n}, b^{+} \in \mathbb{R}^{n}\right\} \\
c(y(b), b)=0 \\
\psi_{j}(b)=\tilde{\psi}_{j}(u(b), y(b)), j=0,1, \ldots, m
\end{matrix}
\label{eq:discretizedmodel}
\end{equation}

\section{Surrogate Model for Optimization Problem}


Solution of optimization problem in equation (\ref{eq:discretizedmodel}) requires multiply solutions of state equation $c(y(b), b)=0$ and calculation of optimization criteria $\tilde{\psi}_{0}$ and constraints $\tilde{\psi}_{j}, j=1, \ldots, m$ of the system for different values of optimization parameters $b$. Complexity of mathematical models (state equation), which describe state and behavior of considered dynamical system requires significant computing resources (CPU time, memory,...) and occasionally puts in question the solving of the optimization problem itself. In order to solve multi-query problems within limited computational cost, there is a need to construct approximation models (also known as surrogates models, meta-models or ROMs). Surrogate model replaces the high-fidelity problem and tends to much lower numerical complexity.

In this paper surrogate models are constructed by first selecting a sampling strategy. Then, $n_s$ sampling points are generated and for each sample point $b^{(i)}$, we solve state equation in equation (\ref{eq:discretizedmodel}) (ODEs) and obtain $n_s$ vectors of solutions (snapshots) $Y_{i}=\left[y\left(t_{1}, b^{(i)}\right)^\top, \ldots, y\left(t_{n_{t}}, b^{(i)}\right)^\top \right]^\top \in \mathbb{R}^{m},  m=n_{y} \times n_{t} $ at different time instances, $t_{0}<t_{1}<t_{2}<\ldots<t_{n_{t}}=T . $ Snapshots vectors $ Y_{i}$ create snapshot matrix $Y=\left[Y_{1}, Y_{2}, \ldots, Y_{n}\right] \in \mathbb{R}^{m \times n_{s}}$.

Next, we construct surrogate model using POD and RBF and calculate the value of functionals $\hat{\psi}_{j}(b)=\tilde{\psi}_{j}(b, \hat{y}), j=0,1, \ldots, m$. Detailed description of POD-RBF procedure is presented in the following paragraphs of this chapter.
The formulation of optimal control problem for surrogate model is to find a vector $\hat{b}^*$ such that:

\begin{equation}
\begin{matrix}

\hat{b}^{*} \in U_{\partial}, \hat{\psi}_{0}\left(\hat{b}^{*}\right)=\displaystyle \min _{b \in U_{\partial}} \hat{\psi}_{0}(b) \\
U_{\partial}=\left\{b: b \in U_{b}, \hat{\psi}_{j}(b)=0, j=1, \ldots, m_{1} ; \hat{\psi}_{j}(b) \leq 0, j=m_{1}+1, \ldots, m\right\} \\
U_{b}=\left\{b: b \in \mathbb{R}^{n}, b^{-} \leq b \leq b^{+}, b^{-} \in \mathbb{R}^{n}, b^{+} \in \mathbb{R}^{n}\right\} \\
\quad \hat{y}=S(b)\\
\hat{\psi}_{j}(b)=\tilde{\psi}_{j}(u(b), \hat{y}), j=0,1, \ldots, m

\end{matrix}
\label{eq:surrogatediscretized}
\end{equation}

Replacing the state equation in (\ref{eq:discretizedmodel}) with surrogate model given in equation (\ref{eq:surrogatediscretized}) is hypothesized to decrease the computational time by a significant amount, because it is free of the complexity of initial problem and involves matrix multiplication that can be accomplished in a much smaller time than solving ordinary differential equations with high fidelity methods. The hypothesis is tested by comparing the accuracy of system responses and time of calculation for both equation (\ref{eq:discretizedmodel}) and equation (\ref{eq:surrogatediscretized}). The detailed procedure for testing of surrogate model is discussed in the next chapter.

\section{Initial Sampling and Method of Snapshots}
\label{sec:initialSampling}

The method of snapshots for POD was developed by Sirovich \cite{Sirovich1987} in 1987. Generally, it comprises of evaluating the model equations for the number of sampling points at various time instances. Each model response is called snapshot and is recorded in a matrix which is collectively called snapshot matrix.

The initial dimension of the problem is equal to the number of snapshots $n_s$ recorded at each time instance $t_i, i=1,...,n_t$. There is no standard method for generating the sampling points. Nevertheless, the choice of sampling method has direct effects on the accuracy of the model and therefore, it is regarded as an autonomous problem. This research briefly explores the initial sampling problem by comparing various classical a-priori methods of sampling. The deeper questions of sampling that relate to the choice of surrogate model, nature of the objective function and analysis are left for the reader to explore from recommended sources such as \cite{Emiliano2013}. 

The main sampling methodology used in the computational procedure is Latin Hypercube Sampling (LHS) and its variant Symmetric Latin Hypercube Sampling (SLHS). LHS is a near-random sampling technique that aims at spreading the sample points evenly across the surface. In statistics, a square grid containing sample positions is a Latin square if and only if there is only one sampling point in each row and each column. A Latin hypercube is the generalization of this concept to an arbitrary number of dimensions, whereby each sample is the only one in each axis-aligned hyperplane containing it. Unlike Random Sampling (RS), which is frequently referred as Monte-Carlo method in finance, LHS uses a stratified sampling techniques that remembers the position of previous sampling point and shuffles the inputs before determining the next sampling points. It has been considered to be more efficient in a large range of conditions and proved to have faster speed and lower sampling error than RS \cite{Lonnie2014}.

\begin{figure}[h]
\includegraphics[width=12cm, height=3cm]{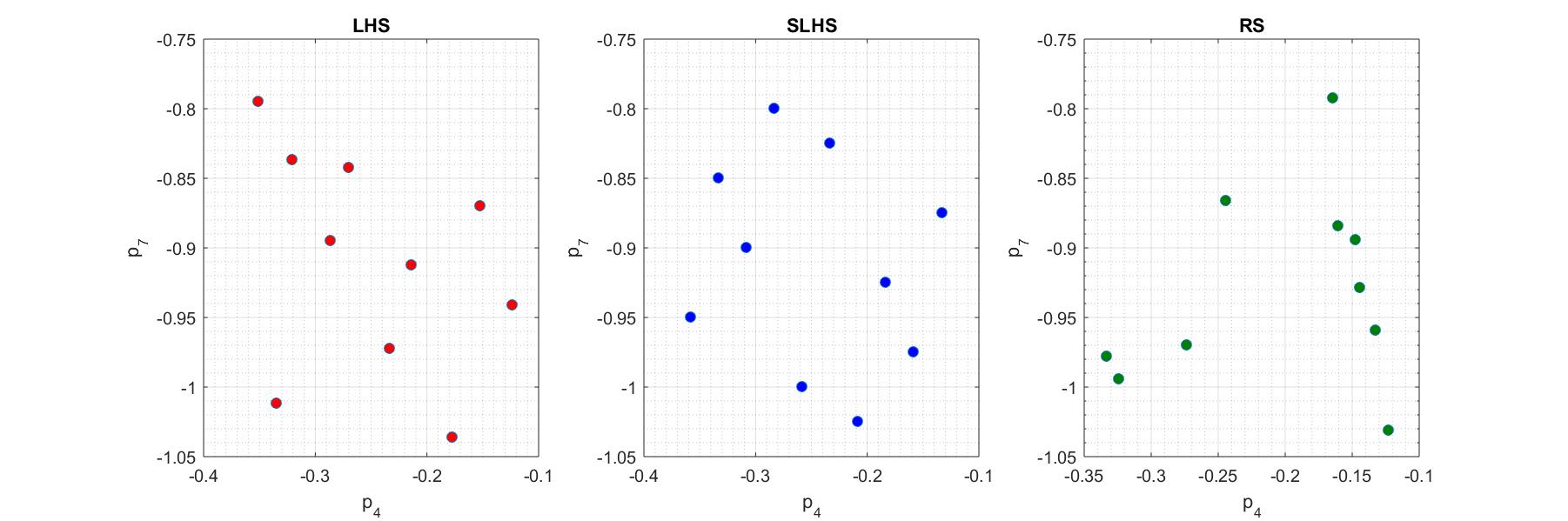}
\caption{Comparison of various sampling techniques.
 \protect \includegraphics[height=0.4cm]{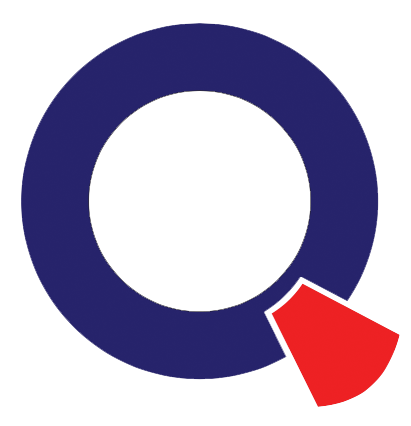} {\color{blue}\href{https://github.com/QuantLet/SurrogateModel}{SurrogateModel}}
}
\label{fig:sampling}
\end{figure}

SLHS was introduced as an extension of LHS that achieves the purpose of optimal design in a relatively more efficient way. It was also established that sometimes, SLHS had higher minimum distance between randomly generated points than LHS. In a nutshell, both LHS and SLHS are hypothesized to perform better than RS. Nevertheless, sampling is performed using all three techniques  in this work to determine which techniques provides optimal sampling of the underlying space and maximizes the system accuracy. A simple sampling distribution of each of the three techniques is illustrated in figure \ref{fig:sampling}.

\vspace*{-10mm}

\section{Approximation}
The overarching goal of POD method is to provide a fit of the desired data by extracting interpolation functions from the information available in the data set. Geometrically, it derives ROMs by projecting the original model onto the reduced space spanned by the POD modes \cite{Emiliano2013}. A simple mathematical formulation of POD technique is laid out in this section which closely follow the references \cite{bujlak2012,chatterjee2000,shcherbatyy2018}.

Suppose that we wish to approximate the response of the system given by output parameters $y \in \mathbb{R}^m$, where $m=n_y \times n_t$, using the set of input parameters $b \subset \mathbb{R}^{n_u}$ over a certain domain $\Omega$. The ROMS approximate the state function y(t) in domain $\Omega$ using linear combination of some basis function $\phi^i\left(x\right)$ such that
\begin{equation}
y\left(t\right)\approx\sum_{i=1}^{M}{a_i.\phi^i\left(t\right)} 
\label{eq:approx}
\end{equation}
where, $a_i$ are unknown amplitudes of the expansions and \textit{t} is the temporal coordinate. The first step in this process would be to find the basis and choice is clearly not unique. Once the basis function is chosen, the amplitudes are determined by a minimization process and the least square error of approximation is calculated. It is ideal to take orthonormal set as the basis with the property 
\begin{equation}
\int_{\Omega}{\phi_{k_1}\left(t\right).\ \phi_{k_2}\left(t\right)dx=\left\{\begin{matrix}1&k_1=k_2\\0&k_1\neq k_2\\\end{matrix}\right.}
\end{equation}
This way, the determination of the amplitudes $a_k$ only depends on function $\phi_k^i(t)$ and not on any other $\phi$. Along with being orthonormal, the basis should approximate the function in best possible way in terms of the least square error. Once found, these special ordered orthogonal functions are called the POMS for the function y(t) and the equation (\ref{eq:approx}) is called the POD of y(t).

In order to determine the number of POMs that should be used in approximation of lower dimensional space, we use the idea that POD inherently orders the basis elements by their relative importance. This is further clarified in the context of SVD in the next section.

\vspace*{-5mm}
\section{Singular Value Decomposition}

There prevails a misconception amongst researchers about distinction between SVD and POD. As opposed to the common understanding, POD and SVD are not strictly the same: the former is a model reduction technique where as the latter is merely a method of calculating the orthogonal basis. Since the theory of SVD is so widespread, this section will only highlight the most general and relevant details of SVD that are helpful in derivation of POMs and POD basis.

In general, SVD is a technique that is used to decompose any real rectangular matrix \textit{Y} into three matrices, \textit{U}, $\Sigma$ and \textit{V}, where \textit{U} and \textit{V} are orthogonal matrices, where $\Sigma$ is a diagonal matrix that contains the singular values $\sigma_i$ of \textit{Y}, sorted in a decreasing order such that $\sigma_1 \geq \sigma_2 \geq...\geq\sigma_d\geq 0$, where \textit{d} is the number of non-zero singular values of \textit{Y}.

The singular values can then be used as a guide to determine the POD basis. If a \textit{k}-dimensional approximation of original surface is required, where \textit{k<d}, the first \textit{k} columns of the matrix \textit{U} serve as the basis $\phi^i, i=1,...,k$. These set of columns, gathered in matrix $\Phi$, form an orthonormal set of basis for our new low-dimensional surface and \textit{k} is referred as the rank.

After collection of basis using SVD, it is easy to calculate the matrix of amplitudes $A_k$. Let $\Sigma_k=[\sigma_1, \sigma_2,..., \sigma_k]$ be the set of \textit{k} largest singular values of our initial matrix \textit{Y}, then, the matrix of amplitudes is given by $Y_k=\Sigma_k A_k, \  A_k=\Sigma_k^\top Y_k$.

Literature on SVD has established that the relative magnitude of each singular value with respect to all the others give a measure of importance of the corresponding eigen-function in representing elements of the input collection. Based on the same idea, a common approach for selection of number of POMs (\textit{k}) is to set a desired error margin $\epsilon_{\text{POD}}$ for the problem under consideration and choose \textit{k} as a minimum integer such that the cumulative energy \textit{E(k)} captured by first \textit{k} singular values (now POMs) is less than 1-$\epsilon_{\text{POD}}$, i.e. 
\begin{equation}
E(k)=\frac{\displaystyle\sum_{i=1}^k \sigma_i^2}{\displaystyle\sum_{i=1}^d \sigma_i^2}\leq 1-\epsilon^2_{\text{POD}}
\label{eq:commulativeenergy}
\end{equation}

\section{Radial Basis Functions}
With the basis vectors and amplitude matrix, using POD discrete theory, low dimensional approximation of our problem has been constructed. However, the formulation is not very useful since our new model can only give the responses of the system for a discrete number of parameter combinations (those that were previously used to generate the snapshot matrix). Since, in many practical applications (for optimization and inverse analysis), even though the values of input parameters may sometime fall in a particular range, they can be any arbitrary combination of numbers between those ranges. That is why, the newly constructed model needs to be approximated in a better way. In this research, POD is coupled with RBF to create low-order parameterization of high-order systems for accurate prediction of system responses. 

RBF is a unique interpolating technique that determines one continuous function that is defined over the whole domain. It is a widely used smoothing and multidimensional approximation technique. For construction of surrogate model using our current basis, a function $f(b)=y$, where $b$ is the vectors of some parameters and y is the output of the system that has to be estimated. Let $Y_k$ be the reduced dimensional matrix calculated by multiplication of basis and amplitudes matrices. It is now easy to apply RBF to reduced dimensional space where system responses are expressed as amplitudes in the matrix $A_k$. Therefore, the surrogate model takes the form $f_a(b)=a$, where $a$ is the vector of amplitudes. Hence,

\begin{equation}
f(b)=y=\Sigma_kA_k=\Sigma_k f_a(b)=\phi f_a(b)\label{fp}
\end{equation}

When RBF is applied for the approximation of $f_a$, $f_a$ is written as linear combination of some basis functions $g_i$ such that
\begin{equation}
f_a(b)=\left[\begin{matrix}
a^i_{1}\\a^i_{2}\\ .\\.\\.\\a^i_{K}\end{matrix}\right]=\left[\begin{matrix}
d_{11}\\d_{21}\\ .\\.\\.\\d_{K_1}
\end{matrix}\right].g_1(b)+
\left[\begin{matrix}
d_{12}\\d_{22}\\ .\\.\\.\\d_{K_2}
\end{matrix}\right].g_2(b)+
...+
\left[\begin{matrix}
d_{1N}\\d_{2N}\\ .\\.\\.\\d_{K_N}
\end{matrix}\right].g_N(b)
=D.g(b)
\end{equation}

Once the basis functions $g_i$ are known, the aim is to solve for the interpolation coefficients that are collectively stored in matrix B. Since we already have the value of amplitudes $A$ from last step, matrix B can be easily obtained by using the equation $B=G^{-1}A$. Finally, using equation (\ref{fp}), our initial space y can be approximated by:
\begin{equation}
y\approx\Phi.D.g(b)=\hat{y}
\label{surrogatemodel}
\end{equation}

In this work, linear and cubic spline RBF are used for analysis, given by:

\begin{equation}
\begin{matrix}
\text{linear spline}:\  g_j(b)=||b-b_j||; \quad \text{cubic spline}:\  g_j(b)=||b-b_j||^3;
\end{matrix}
\end{equation}
 
Since matrix $\Phi$ and D are calculated once for all, one only needs to compute the vector $g(b)$ for any arbitrary combination of parameters to obtain system response.

 Replacing the state equation (\ref{stateequations}) with surrogate model given in equation (\ref{surrogatemodel}) is hypothesized to decrease the computational time by a significant amount, because it is free of the complexity of initial problem and involves matrix multiplication that can be accomplished in a much smaller time than solving ordinary differential equations with high fidelity methods. The hypothesis is tested by comparing the accuracy of system responses and time of calculation for both equation (\ref{stateequations}) and equation (\ref{surrogatemodel}). The detailed procedure for testing of surrogate model is discussed in the next section. 

\vspace*{-5mm}
\section{Error Analysis}
\vspace*{-2mm}
The final step in the analysis of surrogate models is to determine how accurate the low-dimensional surrogate model are in determination of the system responses.

This is done by generating $n_g$ sample points of set of parameters P, using the same sampling technique that had been adapted for generation of training test. It must be noted that the newly generated test points are not same as the one used to train the model and hence, the system responses of these points occur in between nodes and are ideal for checking the accuracy of the models. Moving on, the system responses $Y_g=[y_1,y_2,...,y_{n_g}] \in \mathbb{R}^{m\times n_g}$ are obtained using initial numerical method (that solves entire system), and also $\hat{Y}_g=[\hat{y}_1,\hat{y}_2,...,\hat{y}_{n_g}] \in \mathbb{R}^{m\times n_g}$ are recorded using newly constructed surrogate model. Then, the accuracy and error measures are generally calculated using the following formulas:

\begin{equation}
 R^2=1-\frac{\displaystyle\sum_{1}^{n_g}|y_j-\hat{y}_j|}{\displaystyle\sum_{1}^{n_g}|y_j-\overline{y_j}|}
\label{r2}
\end{equation}

\begin{equation}
\text{MAE}=\frac{1}{n_g}\sum_{1}^{n_g}|y_j-\hat{y}_j|
\label{MAE}
\end{equation}

\begin{equation}
\text{MXAE}=\max_{1\leq j\leq n_g}|y_j-\hat{y}_j|
\label{MXAE}
\end{equation}

\begin{equation}
\text{RMAE}=\max_{1\leq i \leq m}\max_{1\leq j\leq n_g}\frac{|y_{ji}-\hat{y}_{ji}|}{y_{ji}}
\label{RMAE}
\end{equation}

All four measures are put to use at various instances in the thesis, for example, coefficient of determination ($R^2$) in equation (\ref{r2}), Mean Absolute Error (MAE) in equation (\ref{MAE}), Maximum Absolute Error (MXAE) in equation (\ref{MXAE}) are evaluated for various rank approximations of SVD, whereas a tolerance threshold for elative Maximum Absolute Error (RMAE) in equation (\ref{RMAE}) is defined for testing the accuracy of optimization results obtained through original and surrogate models.
\chapter{Algorithm}
\label{chapter-computation}
\setcounter{section}{0}

While understanding of mathematical formulation of POD-RBF procedure presented in Chapter \ref{MathematicalFormulation} is essential, its implementation can be quite technical as it involves high-dimensional matrices, a series of functions, complex loops and iterative processes. The idea of this chapter is to give detailed description of the algorithm that was implemented in MATLAB for this research. The whole computational procedure is divided into three parts for simplicity: experimental design, training phase and testing phase. For each part, a section of the chapter is devoted in which importance of the steps of algorithm are discussed and the intricacies of computational procedure are highlighted. Finally, the iterative nature of algorithm is elaborated in Section \ref{sec:iterativeProcess}.

\section{Experimental Design}
For the construction of surrogate model for a dynamical system, the proper definition of the optimal control problem and planning an appropriate experimental design holds high importance since these conditions are hypothesized to reflect on the accuracy of the model. This pivotal decision relies on choice of fixed and variable parameters, values of constants for fixed parameters, the domain of variable/control parameters, number of initial sampling points, number of time-instances for computation of snapshots, the sampling strategy, the interpolation technique, and minimum error of approximation/ stopping criteria.

Because of the inherent dependence of model on the factors enlisted above, the decision about experimental design has to be made before the construction of surrogate model. In this research, various combination of these factors are accounted for to determine which experimental design results in the highest accuracy while satisfying the time constraints for generation of the snapshots. The error of approximation can be defined for accuracy of system responses or computational time or both.

\section{Offline/Training Phase}
The offline phase (training of the model) entails utilization of sampling techniques to generate data, computation of snapshot matrix of  model solutions using variable order methods for solving of ODE (model of dynamical system), obtainment of proper orthogonal modes via singular value decomposition and estimation of POD expansion coefficients that approximate the POD basis via RBFs.

The next step in the analysis is to determine the appropriate number of POD modes to be used in the surrogate model. For that, the orthogonal basis are found using SVD and the error measures (\ref{r2}), (\ref{MAE}), and (\ref{MXAE}) are used to determine the singular values (rank) whose corresponding eigenvectors will used as POD basis. Next, the amplitudes of approximation $a_i$ are computed using the basis vectors $\phi_i, i=1,...,k$ and stored in amplitude matrix $A_k$.
With this, the dimensionality of this problem cut from $n_s$ to just \textit{k} (rank). Now, to obtain the system response for any arbitrary data point, it is sufficient to simply multiply the reduced basis with corresponding amplitude.

In the final part of offline phase, POD is combined with RBF. The coefficients of RBF interpolation collected in matrix \textit{D} are calculated using our initial data points in $u$ and our final matrix of amplitudes $A_k$ as inputs. With this, the training phase comes to an end. Now, for the computation of system responses, $y \approx \phi.D.g(b)$, surrogate model can be used with only $g(b)$ remaining to be calculated, which depends upon the test points.

\section{Online/Testing Phase}

The last step in construction of low dimensional model is to check the overall error of approximation. It is done by taking the sample points and for each of the data point, first the original response of the function is recorded by solving the ODEs using MATLAB solver \textit{ode15s}, and then the newly developed surrogate model is used to calculate the system response for the same data point. Finally the error measure RMAE (\ref{RMAE}) is calculated for each experimental design and compared to determine which combination meets the required tolerance level.

After deciding the final sampling strategy, number of sampling points, and interpolation technique, the optimization problem is solved using system responses for both original and surrogate models. Then, RMAE is calculated to evaluate the accuracy of surrogate model. If the accuracy level is above the decided threshold, the algorithm enters an iterative process that allows decreasing the width of domain of control parameters. A detailed discussion of iterative process is demonstrated in next section.

\section{Iterative Process}
\label{sec:iterativeProcess}

As stated in the previous section, when the optimization results are obtained using both original and surrogate model, sometimes the desired accuracy of the model is not obtained in the first iteration, despite selecting the best experimental design. This is because the optimal values are usually the corner points and the predictive models in general tend to perform poorly on extreme ends. One of the most effective method to overcome this issue is by the use of adaptive sampling, a method that has been used by many researchers such as \cite{Emiliano2013} with the aim of finding optimal design space points. Despite the effectiveness of the approach, it was not adopted in this work due to limited scope of the research, as previously explained in Section \ref{sec:initialSampling}.

The algorithm used for this research, on the other hand, caters to the aforementioned issue in two ways. Firstly, it trains the initial model with the sampling points from a slightly wider domain than the domain in which the optimization is performed. This way, the corner points are incorporated into the sampling space and surrogate model tends to provide better approximation for the optimal points. Secondly, in order to minimize the error of approximation, the algorithm allows to decrease the width of domain of control parameters at each iteration. By decreasing the size of design space, the sampling points move closer and even if the corner points are not accounted for in the sampling design, the smallest distance between the corner and the neighboring points is lower in smaller domain, hence resulting in better approximation and higher accuracy. If the accuracy is not achieved, the iterative algorithm becomes active: every time the error of approximation is higher than the tolerance level, it shortens the domain, and reconstructs the surrogate model for analysis. The iterative process can be summarized in four steps:

\begin{enumerate}

\item  Initialization: In this step parameters of algorithm are initialized that are required for the iterative process, such as width (the length of domains of control parameters), desired tolerance level and $b^{(0)} =$ initial guess for \textit{b} (the optimization parameters)

\item Setting up the bounds: In this step, upper and lower bounds of domain are defined for each control parameter. It is done by taking $b^{(0)}$, interpolating it and substituting it as the value of control variables in the data structure. Next, the new bounds are created centered at $b^{(0)}$. The width of domain for each subsequent iteration is lower than the previous iteration. The value of $b^{(0)}$ is replaced with optimal value of \textit{b} obtained using surrogate model ($\hat{b}^*$) in the previous iteration. Finally, it is checked if the new bounds are within the bounds that were defined at the beginning of the problem. If not, the algorithm restricts them to exceed the initial bounds. The step 2 of iterative process is depicted for two optimization parameters in figure \ref{IterationAlgorithm}.

\begin{figure}
\begin{center}\includegraphics[width=12cm, height=6cm]{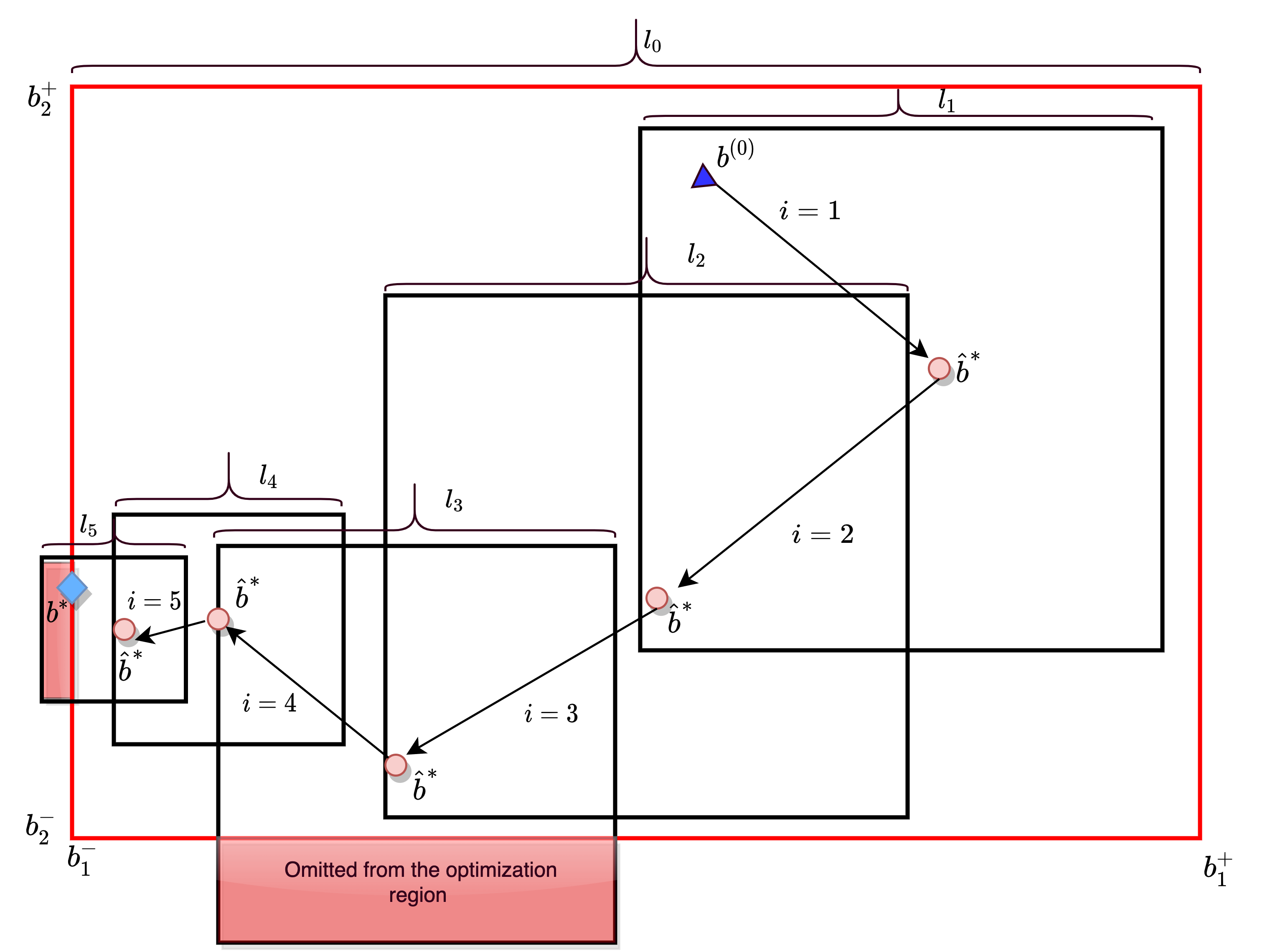}
\caption{Example of iterative algorithm of two optimization parameters $b_1 \text{ and } b_2$ with iterations $i= 1, \ldots, 5$ and recursively decreasing lengths $l_i, i= 1,\ldots,5 $}
\label{IterationAlgorithm}
\end{center}
\end{figure}

\item Optimization: This is the main step of algorithm which was discussed in detail in the second and third section of this chapter. In summary, we make sampling set and snapshots, create surrogate model, solve optimization problem and calculate error. 

\item Updating parameters: This step prepares the parameters for the next iteration in the case when the tolerance level falls below the error of approximation. In general, the algorithm  replaces $b^{(0)}$ with the optimized value of $\hat{b}^*$ from the surrogate response of current iteration, shortens the length by using a predefined multiplier. If the tolerance criteria is met, it stops the iterative process. Else it goes back to step 2. 
\end{enumerate}

\begin{figure}
			\begin{center}				\includegraphics[width=12cm, height=4.5cm]{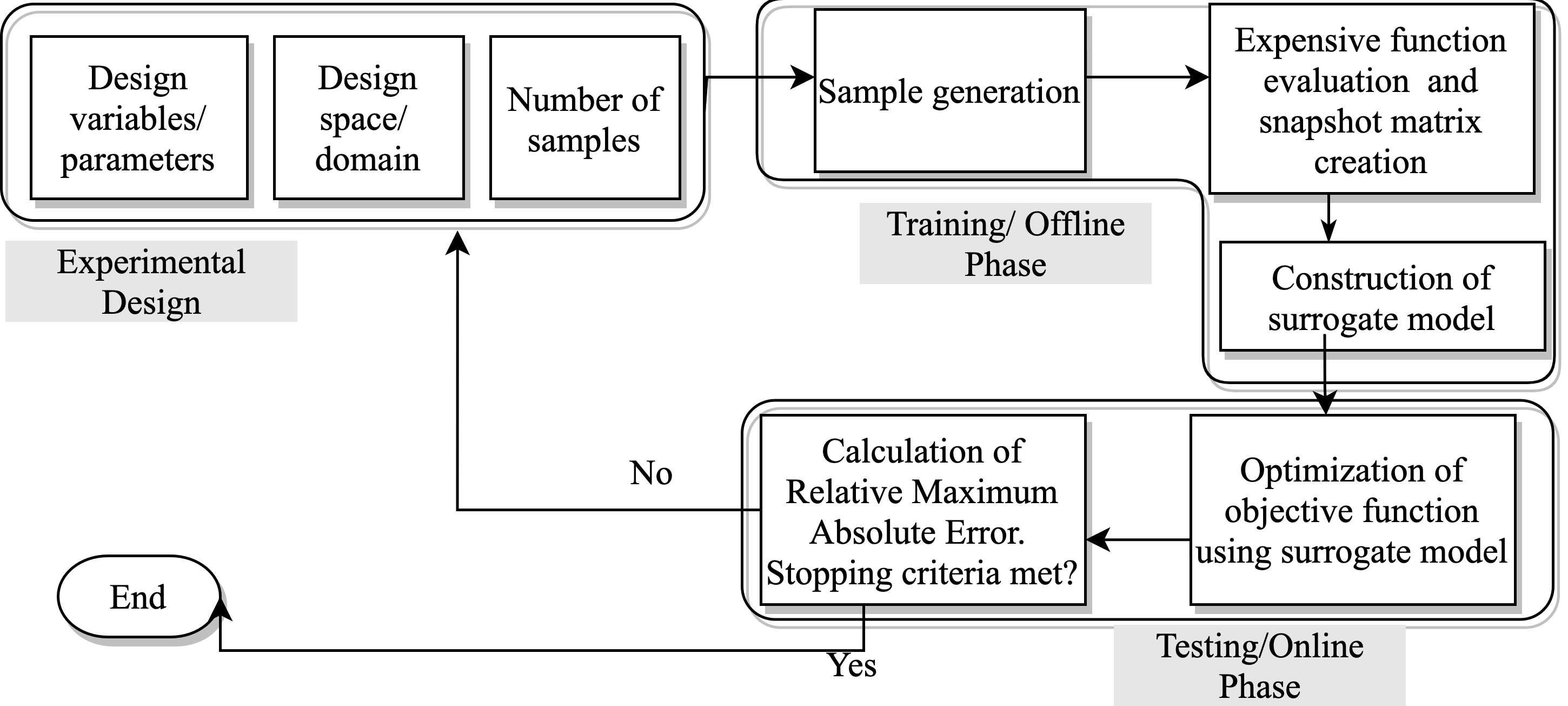}
				\caption{POD-RBF algorithm flowchart}
				\label{fig:flowchart}
			\end{center}
		\end{figure}

The computation procedure discussed throughout this chapter is summarized in flowchart presented in the figure \ref{fig:flowchart}. 
\chapter{Application of POD-RBF Procedure on Dynamical Systems}
\label{chapter-applications}
\setcounter{section}{0}

In this chapter, the POD-RBF procedure is trained and used to construct the surrogate models for real-life dynamical systems and solve associated optimization problems. Two dynamical systems with various complexity are presented, with model 1 being the simple non-linear ODE problem, and model 2 featuring a non-linear system of equations with complex optimization criteria. For each model, a description of the problem is presented and the values of initial parameters used in numerical experiments are defined. Next, the numerical experiments are performed to first decide the combination of sampling technique, interpolation method and sampling points optimal for that model and then the optimization problem is solved to evaluate the accuracy of surrogate responses and the difference in computational time of optimization with original and POD-RBF methods.

As a convention for this chapter, the variables with hat represent the results obtained using surrogate model and without hat stand for the results from original model. The description of common variable names are summarized in table \ref{notations}.

\begin{table}[]
\centering
\begin{tabular}{ll}
\hline
\hline
\textbf{Notation}\quad\quad & \multicolumn{1}{l}{\textbf{Description}} \\ 
\hline
$b^{(0)}$ & Initial value of optimization parameter \\ 
$\hat{b}^*$ & \begin{tabular}[c]{@{}c@{}}Optimal value of optimization parameter, surrogate model\end{tabular} \\ 
$b^*$ & \begin{tabular}[c]{@{}c@{}}Optimal value of optimization parameter, original model\end{tabular} \\ 
$\psi_0(b^{(0)})$ & Value of optimization criteria for $b^{(0)}$, original model\\ 
$\psi_0(\hat{b}^*)$ &  Value of optimization criteria for $\hat{b}^*$, original model \\ 
$\widehat{\psi_0}(\hat{b}^*)$ &  Value of optimization criteria for $\hat{b}^*$, surrogate model \\ 
$\psi_0(b^*)$ &  Value of optimization criteria for $b^*$, original model  \\ 
$\psi_i(b^{(0)})$ &  Value of $i^{th}$ optimization constraint for $b^{(0)}$, original model \\ 
$\psi_i(\hat{b}^*)$ &  Value of $i^{th}$ optimization constraint for $\hat{b}^*$, original model \\ 
$\widehat{\psi_i}(\hat{b}^*)$ &  Value of $i^{th}$ optimization constraint for $\hat{b}^*$, , surrogate model \\
$\psi_i(b^*)$ &  Value of $i^{th}$ optimization constraint for $b^*$, original model  \\ 
\hline
\hline
\end{tabular}
\caption{Details of notations used in preceding analysis}
\label{notations}
\end{table}

\section{Model 1: Science Policy}

\subsection{Description of the Model}
This section features a very interesting application of optimal control theory in economics. The problem is one of the oldest optimal control problem in economics known as science policy and was originally introduced in 1966 by M.D. Intriligator and B.L.R. Smith in their paper "Some Aspects of the Allocation of Scientific Effort between Teaching and Research" \cite{Intriligator1966}. Science policy addresses the important issue of allocation of new scientists between teaching and research staff, in order to maintain the strength of educational processes or alternatively, avoiding any other dangers caused by inappropriate allocation between scientific careers \cite{Intriligator1975}. In order to find the optimal allocation, the optimal control problem was formulated as following:

\begin{equation}
\max _{(u,y) \in U \times Y} \tilde{\psi_0} =\int_{t_0}^\top [0.5y_1(t)+0.5y_2(t)]dt,
\label{model2}
\end{equation}

\begin{equation*}
\begin{matrix}
\text{subject to} \quad \quad \quad\quad\quad\quad \quad \quad \quad\quad\quad\quad \quad \quad \quad \quad\quad\quad\quad  \\
c(y,u)=0 \sim \begin{dcases}
y_1'(t)- u(t)gy_1(t)+\delta y_1(t)=0, t \in [t_0,T]\\
y_2'(t)-(1-u(t))g y_1(t)+\delta y_2(t)=0\\
y_1(t_0)-y_{10}=0, y_2(t_0)-y_{20}=0
\end{dcases}\\
\begin{bmatrix}
  \tilde{\psi_1}  \\
  \tilde{\psi_2}  
\end{bmatrix}= \begin{bmatrix}
  0\\
  0
\end{bmatrix} \sim \begin{dcases}
y_1(T)-y_{1T}=0\quad \quad \quad\\
y_2(T)-y_{2T}=0\quad\quad\quad\quad\quad\quad\quad\quad\quad\quad \quad
\end{dcases} \\
u^-\leq u(t) \leq u^+\quad \quad \quad\quad\quad\quad \quad \quad \quad\quad \quad 
\end{matrix}
\end{equation*}

In this formulation, the state variable $y_1$ and $y_2$ represent the teaching scientists and research scientists respectively at any given time t. The detailed description of all the parameters and their values are summarized in table \ref{model2parameters}. As the control variable $u$ represents the number of new scientists becoming teachers, $(1-u)$ represents the proportion of researchers. Hence, the differential equations determine the rate of change of number of teachers and researchers by subtracting the exiting proportion from the allocated proportion. The upper and lower limit of control function indicate the limits of the science policy in affecting the initial career choices, by government contracts, grants, incentive schemes, etc.

\begin{table}[ht]
\centering
\begin{tabular}{llr}
\hline
\hline
\textbf{\begin{tabular}[c]{@{}c@{}}Parameters \end{tabular}} & \textbf{Definitions} & \textbf{Values} \\ \hline
$u(t_0)$  & \begin{tabular}[c]{@{}c@{}} Proportion of new scientists becoming teachers at initial time \end{tabular}& 0.5 \\ 
$g$ & \begin{tabular}[c]{@{}c@{}} Number of scientists annually produced by one scientist \end{tabular}  & 0.14 \\ 
$\delta$ & \begin{tabular}[c]{@{}c@{}} Rate of exit of scientists due to death, retirement or transfer \end{tabular} & 0.02 \\ 
$y_{10}$ & Number of initial scientists working as teachers  & 100 \\ 
$y_{20}$ & Number of initial scientists working as researchers    & 80 \\ 
$T$ & Final time for the analysis in this policy  &  15\\ 
$y_{1T}$ & Number of final scientists working as teachers  & 200 \\ 
$y_{2T}$ & Number of final scientists working as researchers    & 240 \\ 
$u^-$ & Lower limit of control function& 0.1 \\ 
$u^+$ & Upper limit of control function & 0.6 \\ \hline
\hline
\end{tabular}
\caption{Description of parameters for Model 1}
\label{model2parameters}
\end{table}

The problem is the one of choosing a trajectory for the allocation of $u(t)$ such that the welfare is maximized, given by the objective function in equation (\ref{model2}). The terminal part $g_1(.,.)$ of welfare is not accounted for in the objective function, but the state constraints are added to compensate for it in the form of $y_1(T)-y_{1T}=0$ and $y_2(T)-y_{2T}=0$. The optimization process is focused at maximizing the intermediate value $g_2(.,.,.)$ of welfare. The welfare function is thought to be additive of individual utilities along the lines of utilitarian approach. The utilities are set as a linear function, with an assumption that the teachers and researchers are perfect substitutes, and the allocation of any scientist to one career will lead him to abandon the other career completely. The assumption, even though unrealistic, is granted for simplicity and can be complicated at the later stages.

\subsection{Simulation}
\begin{figure}
\begin{center}
\includegraphics[width=12cm, height=4cm]{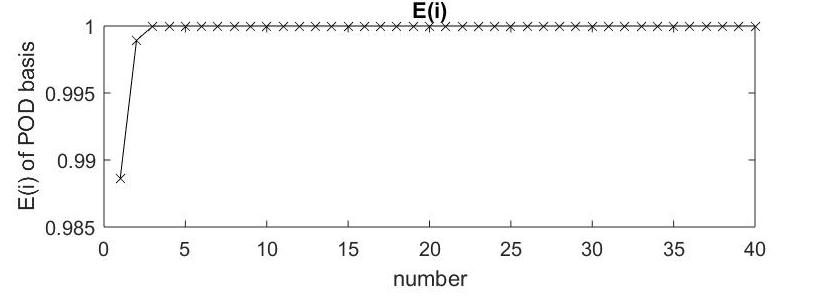}
\vspace*{-10mm}
\caption{Cumulative energy plot to determine singular values for Model 1. \protect \includegraphics[height=0.4cm]{figures/qletlogo_tr.png} {\color{blue}\href{https://github.com/QuantLet/SurrogateModel/tree/main/SurrogateModel\_SciencePolicy}{SurrogateModel\_SciencePolicy}}}
\label{fig:model2singularvalues}
\end{center}
\end{figure}

This system of equation is solved for $n_s=40,60,80$ training points, generated with LHS, SLHS and RS to create the snapshot matrix. The desired tolerance level is $\epsilon_{\text{POD}}=0.01$. The plots of singular values of depicted similar pattern for all the experimental designs. The singular value plot for one specific example, SLHS and $n_s=40$ is presented figure \ref{fig:model2singularvalues} and shows that the first 4 singular values explain almost 100\% variance. Given the criterion in equation (\ref{eq:commulativeenergy}), we choose the rank of $k=4$.

\begin{figure}[h]
\begin{center}
\includegraphics[width=12cm, height=6cm]{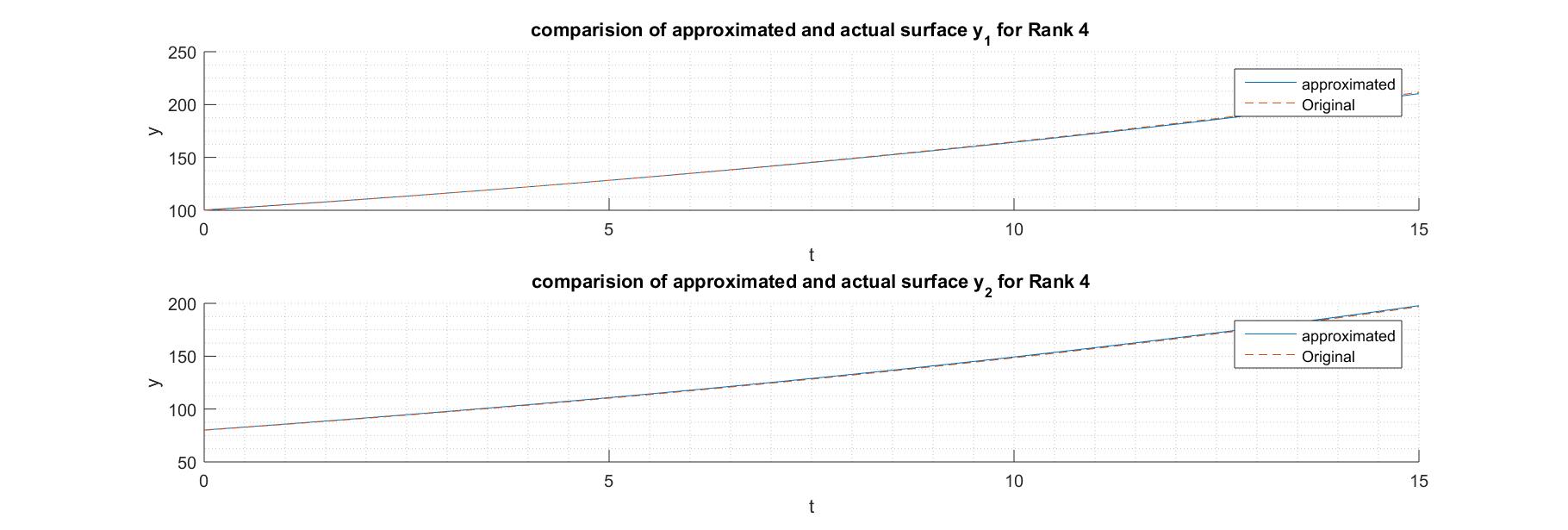}
\vspace*{-15mm}
\caption{Actual surface vs approximated surface for Model 1.
\protect \includegraphics[height=0.4cm]{figures/qletlogo_tr.png} {\color{blue}\href{https://github.com/QuantLet/SurrogateModel/tree/main/SurrogateModel\_SciencePolicy}{SurrogateModel\_SciencePolicy}}}
\label{fig:model2actualvsapprox}
\end{center}
\end{figure}

The surrogate model was constructed for each of the variant with this rank and the RMAE are reported in table \ref{tab:model2RMAE}. The table shows that the lowest RMAE was obtained for LHS, followed by SLHS and the RS. As the theory suggests, RMAE is observed to decrease with increasing number of sampling points with an exception of cubic spline in random sampling. The anomalous behavior of RS can be associated with its randomness, which sometimes generates the sampling points which belong to only one region of the surface, leading to higher variance in the model and higher error of approximation, even with increasing number of training points. Another trend that can be consistently observed is that the linear spline RBF tend to perform better than the cubic spline in this model. Overall, the best experimental design for this model is to use a combination of LHS with linear spline RBF and $n_s=80$. The surrogate model approximation for the initial control value $u=0.5$ and the original system response are plotted in figure \ref{fig:model2actualvsapprox} and show that the approximated responses are very close to the actual responses.

\begin{table}[ht]
\centering
\begin{tabular}{llllr}
\hline
\hline
\textbf{Sampling\quad\quad } & \textbf{Interpolation\quad\quad } & \textbf{$n_s=40$} \quad\quad   & \textbf{$n_s=60$}  \quad\quad   & \textbf{$n_s=80$}\quad\quad\\ \hline
\multirow{2}{*}{LHS} & Linear & 0.02034 & 0.00293 & 0.00150 \\ 
 & Cubic & 0.05316 & 0.00647 & 0.00641 \\ 
\multirow{2}{*}{SLHS} & Linear & 0.03825 & 0.00679 & 0.00437 \\ 
 & Cubic & 0.05175 & 0.00897 & 0.00861 \\ 
\multirow{2}{*}{RS} & Linear & 0.01525 & 0.02410 & 0.02792 \\  
 & Cubic & 0.16457 & 0.26597 & 12.91601\\ 
\hline
\hline

\end{tabular}
\caption{RMAE for various experimental designs of Model 1}
\label{tab:model2RMAE}
\end{table}

\vspace*{-10mm}

\subsection{Optimization}

\begin{figure}[h]
\begin{center}
\includegraphics[width=12cm, height=5cm]{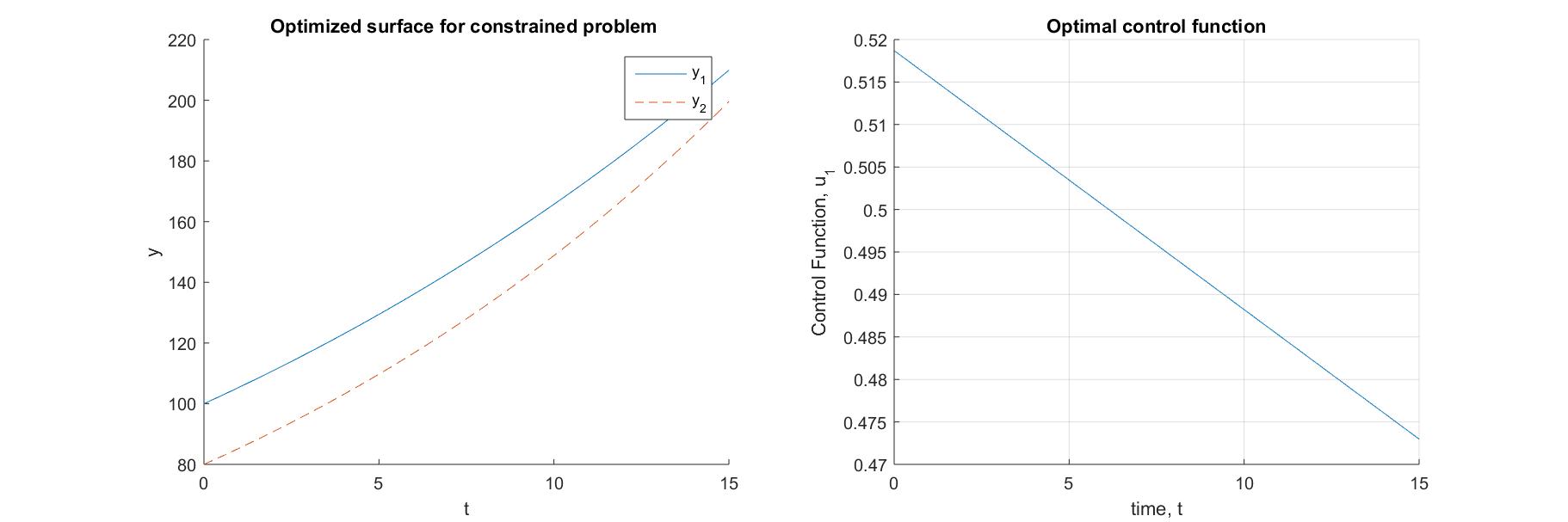}
\vspace*{-10mm}
\caption{Optimal surface and control functions for Model 1. 
\protect \includegraphics[height=0.4cm]{figures/qletlogo_tr.png} {\color{blue}\href{https://github.com/QuantLet/SurrogateModel/tree/main/SurrogateModel\_SciencePolicy}{SurrogateModel\_SciencePolicy}}}
\label{fig:model2optimization}
\end{center}
\end{figure}

For the final step of analysis, the surrogate model was constructed with 40 training points, LHS, and linear spline RBF. Here $n_s=40$ was used because given the simplicity of the problem, the accuracy required for optimization can be achieved by small number of training points. The optimization problem is solved with two optimization parameters for control function using both original and surrogate model. The results of optimization are given in table \ref{tab:model2optimization}. The problem started with equal number of scientists allocated in both careers, with the initial value of state constraint $\psi_1(b^{(0)})=[11.8001;43.0163]$ representing that the number of teachers and researchers allocated at initial time were 11 and 43 units short of $y_{1T}$ and $y_{2T}$ respectively. The solution to the problem allocates around 52\% of new scientists to teaching at the beginning of the time. This proportion decreases as the time passes with around 47\% scientists allocated as teaching staff at the end of time (see figure \ref{fig:model2optimization}(b)). The optimal surface in \ref{fig:model2optimization}(a)) shows that the number of teaching staff was allocated to be higher than the number of researchers until the end time. The surrogate model gave consistent results, with error of approximation (the relative error of $\psi_0(\hat{b}^*)$ and $\widehat{\psi_0}(\hat{b}^*)$) as low as 0.005 in the first iteration. 

Even though the optimization using surrogate model was slightly quicker than the original model, the time taken for construction of surrogate model was higher. Hence, despite of highly accurate system responses through surrogate model, substituting original model with POD-RBF model might not be useful, as the time taken for optimization by surrogate model (training + optimization) was much longer than the original model. This example give us insight into why surrogate modelling was avoided into applications earlier: the simple nature of optimization models for some applications do not require high computational resources, while the construction of surrogate models is much more computationally expensive and may not be desirable.

\begin{table}
\centering
\begin{tabular}{llll}
\hline
\hline
\textbf{Field} & \textbf{Value} & \textbf{Field} & \textbf{Value} \\
\hline 
$b^{(0)}$ & {[}0.5000 0.5000{]} & Bounds & {[}0.1000,0.6000{]} \\ 
$b^*$ & {[}0.6000,0.3461{]} & $\hat{b}^*$ & {[}0.5187,0.4730{]} \\
$\psi_0(b^{(0)})$ & 210.6500 &  $\psi_0(\hat{b}^*)$ & 209.7600\\
$\psi_0(b^*)$ &212.8400& $\widehat{\psi_0}(\hat{b}^*)$ &210.9900 \\
$\psi_1(b^{(0)})$ & ${[}11.8001,43.0163{]}^\top$ & $\psi_1(\hat{b}^*)$ & ${[}0.0003,0.0014{]}^\top$ \\ 
$\psi_1(b^{*})$ & ${[}0.000,0.000{]}^\top$  & $\widehat{\psi_1}(\hat{b}^*)$ & ${[}0.0000,0.0023{]}^\top$ \\ 
$\text{Time}_{orig}$ & 2.8109 sec & $\text{Time}_{surr}$ &2.3694 sec \\ 
$\text{Time}_{cnstr}$& 37.8406 sec &$\epsilon$ & 0.0058 \\ \hline
\hline
\end{tabular}
\caption{Optimization results of Model 1}
\label{tab:model2optimization}
\end{table}

\section{Model 2: Population Dynamics}
\subsection{Description of the Model}
In this section, a more complex application of optimal control theory is presented with a general model of non-linear system of ODEs defined by:

\begin{equation} 
c(y, u)=0 \sim\left\{\begin{array}{l}\left\{\begin{array}{l}y_{1}^{\prime}-p_{1} y_{1}-p_{2} y_{2}^{2}-u_{1} y_{1} F\left(y_{1}, t\right) y_{2}=0, \\ y_{2}^{\prime}-p_{3} y_{2}-p_{4} y_{2}^{2}-u_{1} u_{2} y_{1} F\left(y_{1}, t\right) y_{2}=0,\end{array} t \in \Omega_{t}=\left(t_{0}, T\right]\right. \\ y_{1}\left(t_{0}\right)-y_{10}=0 \\ y_{1}\left(t_{0}\right)-y_{20}=0 \\ F\left(y_{1}, t\right)=1-e^{-p_{5} y_{1}}\end{array}\right.
\label{model5}
\end{equation}


These type of dynamical problems are usually observed in population dynamics in biology, ecology and environmental economics. These problems are variation of prey-predator model presented by Lotka-Volterra. This section aims at generalizing the approach of POD-RBF on these non-linear models without providing specific details of the model parameters of the optimization problem.

The optimization problem considered here consists of finding a value of control function $u^{*}=\left[u^*_1, u^*_{2}\right]$ that minimizes the distance between $y_{1}$ and its desirable value $y_{1d}$ Value on control function is restricted by dual pointwise constraints and value $y_{2}$ do not exceed maximum value $y_{2 d} .$ The optimization problem can be formulated in the following manner: find $u^{*}$ that minimize optimization criterion

\begin{equation}
\psi_{0}\left(u^{*}\right)=\min _{u} \int_{t_{0}}^{T}\left(y_{1}(t, u)-y_{1 d}\right)^{2} dt
\label{eq:model2opt}    
\end{equation}
subject to state equation (\ref{model5}), box constraints on the control
\begin{equation}
U=\left\{u: u^{-}(t) \leq u(t) \leq u^{+}(t)\right\}
\label{eq:boxconstraintmodel2}
\end{equation}
and pointwise constraint on state

\begin{equation}
y_{2}(t) \leq y_{2}^{+}
\label{eq:pointwisesc}
\end{equation}

The pointwise state constraint (\ref{eq:pointwisesc}) is transformed into an equivalent equality constraint of the integral type
\begin{equation}
\psi_{1}(u)=\tilde{\psi}_{1}(u, y(u))=\int_{t_{0}}^{T}\left(\left|y_{2}(t, u)-y_{2 d}\right|+y_{2}(t, u)-y_{2 d}\right)^{2} dt
\label{eq:model2constraints}
\end{equation}
Taking into account equations(\ref{eq:model2opt}-\ref{eq:model2constraints}) the optimization problem can be written in a reduced form as follows:
\begin{equation}
    \begin{aligned}
\psi_{0}\left(u^{*}\right)=& \min _{u \in U_{\partial}} \int_{t_{0}}^{T}\left(y_{1}(t, u)-y_{1 d}\right)^{2} d t \\
U_{\partial u}=\{u&\left.: u \in U ; \psi_{1}(u)=\tilde{\psi}_{j}(u, y(u))=0\right\} \\
c(y(u), u) &=0
\end{aligned}
\label{model5optimization}
\end{equation}
\vspace*{-15mm}

\subsection{Simulation}

For numerical experiments we select the following values of the input parameters:
$\left[p_{1}, p_{2}, p_{3}, p_{4}, p_{5}\right]= [0.734,0.175,-0.500,-0.246,0.635],\left[t_{0}, T\right]=[0,10]$, 
$n_{u}=2, u^{-}=\left[u_{1}^{-}, u_{2}^{-}\right]=[-0.5500,-1.0370], u^{+}=\left[u_{1}^{+}, u_{2}^{+}\right]=[-0.300,-0.7870]$,
 $y_{1 d}=5, y_{2}^{+}=6$.
The control functions $u_{1}(t), u_{2}(t)$ on the interval $\left[t_{0}, T\right]$ are approximated by linear functions. Thus, the vector of optimization parameters $b$ consist of four components:
$b=\left[b_{1}^{(1)}, b_{2}^{(1)}, b_{1}^{(2)}, b_{2}^{(2)}\right]^{T}=\left[b_{1}, b_{2}, b_{3}, b_{4}\right]^{T}$.

For numerical simulations, LHS, SLHS and RS are used to define the sampling matrix with $n_s=40, 60 \ \text{and} \ 80$. Also, RBF interpolation-linear spline and cubic spline is used for comparison of results. The solution $y=[y_1,y_2]$ where $n_y=2$ was then computed for time instances, $t_i$ with $t_0 < t_i<t_{n_t}$,  $n_t=100$ equally spaced instances of t, and $n_s$ sampling points, and then system responses were collected to generate the snapshot matrix. The error of approximation was fixed $\epsilon_{\text{POD}}=0.01$. 

\begin{figure}[]
\begin{center}
\includegraphics[width=12cm, height=3cm]{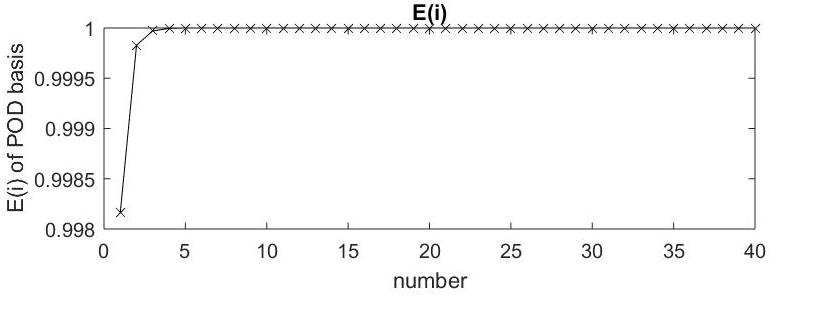}
\vspace*{-10mm}
\caption{Cumulative energy plot to determine singular values for Model 2. 
\protect \includegraphics[height=0.4cm]{figures/qletlogo_tr.png} {\color{blue}\href{https://github.com/QuantLet/SurrogateModel/tree/main/SurrogateModel\_PopulationDynamics}{SurrogateModel\_PopulationDynamics}}}
\label{fig:model5singularvalues}
\end{center}
\end{figure}

Next, the POD-RBF approach is applied to this model to first determine the dimension of POD basis through SVD using cumulative energy method (it is done for all experimental designs) and it is concluded that 3 singular values should be considered as the rank of the POD basis as shown by the figure of singular values in figure \ref{fig:model5singularvalues}. It can be clearly noticed that the magnitude of all the singular values is very small compared to first singular value; the relative commutative energy E(i) of first singular value is more than 99\%. This shows that that the responses of the system are fully correlated. Hence, rank 3 approximation is very accurate and adding more vectors (by increasing rank) in approximation further increases the precision.

\begin{figure}[]
\begin{center}
\includegraphics[width=12cm, height=6cm]{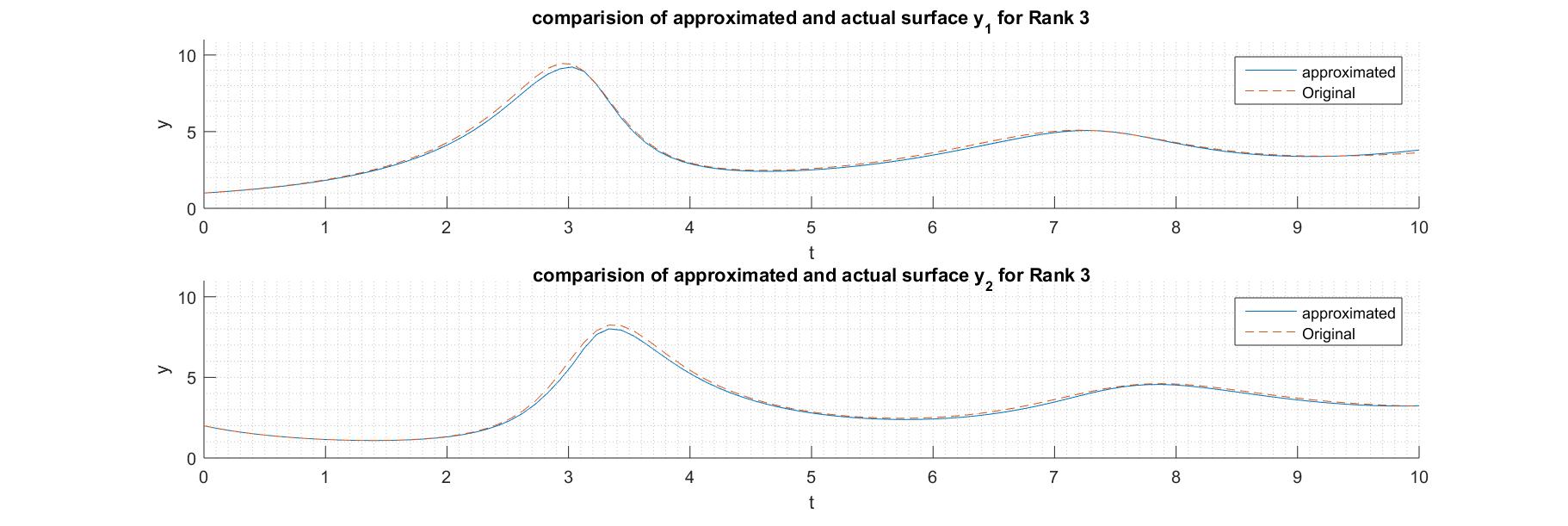}
\vspace*{-10mm}
\caption{Actual vs approximated surface of Model 2. 
\protect \includegraphics[height=0.4cm]{figures/qletlogo_tr.png} {\color{blue}\href{https://github.com/QuantLet/SurrogateModel/tree/main/SurrogateModel\_PopulationDynamics}{SurrogateModel\_PopulationDynamics}}}
\label{fig:p_check}
\end{center}
\end{figure}

Having chosen $k=3$, the numerical simulations are performed for model (\ref{model5}). For testing of the model, $n_g=10$ points were used to calculate the RMAE for each combination. Table \ref{tab:model5RMAE} exhibits that among all the surrogate models that were trained using different number of sample points, different sampling techniques and RBF interpolations, the cubic spline RBF showed the lowest error for both LHS and SLHS in general, with a few exceptions. Also, as expected, the error of approximation shows a decreasing pattern as the number of sample points increase from 60 to 80, except in RS when the RMAE follows no particular trend. The least RMAE was obtained for the model trained on 80 data points from SLHS for cubic spline RBF. For one of such sample point $b=[-0.425,-0.425,-0.912,-0.912]$, the POD-RBF responses were obtained for $n_s=40$ and the original and approximated $y_1$ and $y_2$ were plotted as shown in figure \ref{fig:p_check}. For this point, all POD-RBF gave relative maximum absolute error less than 1\% as desired.

\begin{table}[ht]
\centering
\begin{tabular}{lllll}
\hline
\hline
\textbf{\begin{tabular}[c]{@{}c@{}}Sampling \quad\quad\end{tabular}} & \textbf{\begin{tabular}[c]{@{}c@{}}Interpolation\quad\quad \end{tabular}} & \textbf{$n_s=40$}\quad\quad    & \textbf{$n_s=60$}\quad\quad     & \textbf{$n_s=80$} \quad\quad    \\ \hline
\multirow{2}{*}{LHS} & Linear & 0.45112 & 0.32948 & 0.18871 \\ 
 & Cubic & 0.28229 & 0.24010 & 0.15794 \\ 
\multirow{2}{*}{SLHS} & Linear & 0.26162 & 0.19198 & 0.19204 \\ 
 & Cubic & 0.23986 & 0.18685 & 0.15376 \\ 
\multirow{2}{*}{RS} & Linear & 0.59500 & 0.55080 & 0.86405 \\  
 & Cubic & 0.92109 & 0.15595 & 0.19902\\ 
\hline
\hline
\end{tabular}
\caption{RMAE for various experimental designs of Model 2}  
\label{tab:model5RMAE}
\end{table}

\vspace*{-15mm}

\subsection{Optimization}

In previous subsection, the best results were obtained for $n_s=80$ with SLHS and cubic spline RBF. That experimental design is used to solve the optimization problem (\ref{model5optimization}) and the results are summarized in table \ref{tab:model5optimization}. For simplicity, the number of optimization parameters for each control variable are taken to be 2. We could, however, allows specification of different number of optimization parameters for each control variable. 
The optimization results of this model apparently highlight the efficiency of surrogate modeling. As the table \ref{tab:model5optimization} reports, the tolerance level was met in the first iteration, with error between approximated and actual responses being less than 0.01 in first iteration. Hence, the desired accuracy was achieved and no further iterations were required. Also, the optimization criteria obtained using surrogate model $\widehat{\psi_0}(\hat{b}^*)= 43.5647$ is very close to $\psi_0(b^*)= 43.3287$. Moreover, since results of optimization problem were obtained within one iteration, the construction time of surrogate model can be considered once for all. Therefore, the total computational time for optimization through surrogate model of 6.6 seconds + 15.35 seconds is less than 23.40 seconds taken by original problem. Relatively, the surrogate method was four times faster than the original method in solving optimization problem. In a nutshell, for this highly non-linear model, surrogate model gave highly accurate and computationally efficient result of the optimization problem.  
\FloatBarrier

\begin{table}
\centering
\begin{tabular}{llll}
\hline
\hline
\textbf{Field} & \textbf{Value} & \textbf{Field} & \textbf{Value} \\
\hline 
$b^{(0)}$ & {[}-0.4250,-0.4250,\quad\quad\quad & Bounds & {[}-0.5500, -0.300{]};\\ 
& -0.9120,-0.9120{]}& &{[}-1.0370,-0.7870{]} \\
$b^*$ &[-0.5006,-0.3250, & $\hat{b}^*$ & [-0.4922,-0.3334,\\
&-1.0120,-1.0120]&&-1.0120,-1.0120] \\
$\psi_0(b^{(0)})$ & 55.2817 & $\psi_0(\hat{b}^*)$ & 43.9127\\ 
$\psi_0(b^{*})$ & 43.3287 &$\widehat{\psi_0}(\hat{b}^*)$ &  43.5647\\
$\psi_1(b^{(0)})$ & 22.9396 & $\psi_1(\hat{b}^*)$ & 0.0162  \\ 
 $\psi_1(b^{*})$ & 0.0000& $\widehat{\psi_1}(\hat{b}^*)$ & 0.0000 \\ 
$\text{Time}_{orig}$ & 23.3983 sec & $\text{Time}_{surr}$ & 6.6241 sec \\ 
$\text{Time}_{cnstr}$ & 15.3470 sec & $\epsilon$ & 0.0081\\
\hline
\hline
\end{tabular}
\caption{Optimization results of Model 2}
\label{tab:model5optimization}
\end{table} 

\vspace*{-15mm}
\chapter{Conclusions}
\setcounter{secnumdepth}{0}
This research employed Proper Orthogonal Decomposition (POD), a surrogate modeling technique integrated in optimization framework for dimension reduction by extracting hidden structures from high dimensional data and projecting them on lower dimensional space. In the first instance, POD was coupled with various Radial Basis Functions (RBF)--- a smoothing technique--- and the computational procedure was hypothesized to provide compact, accurate and computationally efficient solution of optimal control problems. The surrogate models using POD-RBF were constructed. The computational procedure of surrogate model was divided into problem setup and training/testing phase for effective implementation of the reduced order modeling techniques. Furthermore, an iterative algorithm was introduced methodically to achieve more accurate results.

The algorithm and computational procedure was implemented on two real-life optimal control problems that were taken directly from literature sources. It was demonstrated that the dimensionality of high order models in the form of ODEs of dynamical systems could be reduced substantially to as low as 3 with relative maximum absolute error less than 0.01 between original and approximated system responses. Hence approximated surrogate model gave a good alternative method of solution of ODEs with low CPU intensity. The simulation part of PDF-RBF procedure was carried out by varying the number of sample points, sampling strategy, and RBF interpolation types in the training phase. The results showed that the approximation was more precise if the model was trained on higher number of sample points. Also, the interpolated surrogate model constructed using cubic-spline RBF led to better results in the complex model than its liner counterpart. Furthermore, LHS and SLHS both led to better approximations than RS which is in coherence with the theory.

In solution of optimization problems, the system responses obtained by surrogate model invariably gave accurate results with improved computational time. As a whole, both the models agreed with the hypothesis of this work that surrogate models can increase the computational efficiency in solution of dynamical systems while maintaining the accuracy of system responses. However, the computational performance is subject to the available computational resources and the numerical simulation might be much faster in a high-performance computer, compensating for the time used in iterative process of POD-RBF algorithm.

\section{Limitations and Future Work} 

ROMs are usually thought of as computationally inexpensive mathematical representations that offer the potential for near real-time analysis. The hypothesis of this research was based on the same notion. However, while analyzing the performance POD-RBF procedure on non-linear dynamical systems in the last chapter of this thesis, it was brought into consideration that the even though the optimization process itself was faster with surrogate responses, their construction was sometimes computationally expensive as it involved accumulating a large number of system responses to input parameters. It is also noteworthy that sometimes ROMs lack robustness with respect to parameter changes. These limitations were considered and elaborated throughout the analysis and the scope of extension of this research was discussed alongside.

In future, the performance of surrogate models can be evaluated on more complicated models consisting of highly non-linear ordinary and partial differential equations. Also, other sampling techniques which allow inclusion of corner and optimization points in the training set, methods of obtaining POMs, and interpolation methods can be explored as an extension of this work. Furthermore, the computational time of each of the model can be calculated with more efficient machines in homogeneous computer environment to get near-exact insight into the performance of surrogate models. 

\renewcommand\bibname{REFERENCES}
 \addcontentsline{toc}{chapter}{References}  
\printbibliography
\pagebreak
\end{document}